# A BAYES METHOD FOR A MONOTONE HAZARD RATE VIA S-PATHS[1]

By Man-Wai Ho

*National University of Singapore*

A class of random hazard rates, which is defined as a mixture of an indicator kernel convolved with a completely random measure, is of interest. We provide an explicit characterization of the posterior distribution of this mixture hazard rate model via a finite mixture of **S**-paths. A closed and tractable Bayes estimator for the hazard rate is derived to be a finite sum over **S**-paths. The path characterization or the estimator is proved to be a Rao–Blackwellization of an existing partition characterization or partition-sum estimator. This accentuates the importance of **S**-paths in Bayesian modeling of monotone hazard rates. An efficient Markov chain Monte Carlo (MCMC) method is proposed to approximate this class of estimates. It is shown that **S**-path characterization also exists in modeling with covariates by a proportional hazard model, and the proposed algorithm again applies. Numerical results of the method are given to demonstrate its practicality and effectiveness.

**1. Introduction.** In reliability theory and survival analysis a hazard rate $\lambda(t)$ is interpreted as the propensity of failure of a system (or an item) in the instant future given that it has survived until time $t$. In general, the function has a wide variety of shapes. The simplest case of a constant hazard rate corresponds to an exponential lifetime distribution for the system. Cases of increasing or decreasing hazard rate, broadly speaking, correspond to lifetime distributions that are of a lighter or heavier tail, respectively, compared to an exponential distribution. There is a substantial amount of literature about estimation of monotone hazard rates from a frequentist viewpoint. They include, for example, the pioneering work of Grenander

Received September 2004; revised May 2005.

[1]Supported in part by National University of Singapore Research Grant R-155-000-047-112 and Hong Kong RGC Competitive Earmarked Research Grant HKUST6159/02P.

[17] and Prakasa Rao [47], extensions of their work to different censoring schemes by Padgett and Wei [46] and Mykytyn and Santner [45], a constrained spline smoothing technique by Villalobos and Wahba [49], work of Lo and Phadia [42] and Huang and Wellner [26] based on the least concave/greatest convex minorants, and a kernel-based method by Hall, Huang, Gifford and Gijbels [19].

A Bayesian nonparametric approach to this important problem is to use the fact that a monotone hazard rate on the half line $\mathcal{R} = (0, \infty)$ may be written in the form

$$\lambda(t|\mu) = \int_{\mathcal{R}} \mathbb{I}(t < u) \mu(du), \tag{1}$$

where $\mathbb{I}(A)$ is the indicator function of a set $A$ and $\mu$ is modeled as a random process. Drăgichi and Ramamoorthi [13] establish the strong and weak consistency of the posterior distribution of these hazard rates for various choices of $\mu$. The consistency of the Bayes estimate follows as a consequence. This is important as it shows that such models yield viable estimators and, hence, alternatives to the approaches mentioned earlier. This approach was first utilized by Dykstra and Laud [14], wherein $\mu$ is modeled as an extended or weighted gamma process. Lo and Weng [41], specifying a weighted gamma process [38] for $\mu$ and replacing the indicator function in (1) by more general kernels $k(t|u)$, provide explicit posterior analysis for hazard rates with more general shapes, that is,

$$\lambda(t|\mu) = \int_{\mathcal{R}} k(t|u) \mu(du). \tag{2}$$

Their analysis, paralleling that of Lo [39], shows that, for general kernels, the posterior distribution can be characterized in terms of random partitions, say, $\mathbf{p} = \{C_1, \ldots, C_{n(\mathbf{p})}\}$, of the integers $\{1, \ldots, n\}$ related to the Chinese restaurant process. Here for exchangeable observations $X_1, \ldots, X_n$, $C_j = \{i : X_i = X_j^*\}$, for $n(\mathbf{p}) \leq n$ unique values $X_1^*, \ldots, X_{n(\mathbf{p})}^*$. This is formally obtained from a posterior distribution which is a mixture of the Pólya urn distribution [7]. See [30] for a recent discussion of these models relative to the Dirichlet and gamma processes in a semi-parametric context. It is now recognized in Bayesian statistics and spatial statistics that other model specifications for $\mu$ may be of interest. Such generalizations of (2), also known as Lévy moving averages, are discussed in, for instance, [50]. James ([32], Section 4), in analogy to Lo and Weng [41], provides a partition-based representation of the posterior distribution for the case where $\mu$ is a completely random measure [36, 37] or related model.

While models in (1) are special cases of (2) and, hence, they have a posterior distribution expressible in terms of random partitions, Dykstra and Laud [14] and Lo and Weng [41] showed that these models actually have a considerably less complex representation in terms of what are called **S**-paths. Quite specifically, **S**-paths are combinatorial structures which, relative to **p**, only contain information about the maximal element and size of each cell $C_j$. This phenomenon is discussed in more detail in the case of monotone densities in [10]. In this work we note the fact that the occurrence of tractable **S**-paths for monotone hazard rates is due to the nice features of the indicator kernel. Hence, using this fact, we are able to refine the partition-based results of James [32] to show that all such monotone hazard rates have **S**-path structures. This represents the first explicit characterization of this type. The main attractive feature is that the space of **S**-paths is considerably smaller than the space of partitions (see Table 2 in the Appendix). Hence, it has been recognized that if one could efficiently sample **S**-paths in this context, this would lead to more parsimonious methods for inference. However, it turns out that the design and implementation of efficient numerical methods utilizing **S**-paths are not that obvious. Section 3 presents an efficient MCMC method for sampling directly **S**-paths induced by monotone hazard rates. We shall also extend this to the semi-parametric setting of a proportional hazard model with covariates in Section 6.

**2. A posterior distribution of a monotone hazard rate model via S-paths.** This section concerns Bayes estimation of a decreasing hazard rate on the half line $\mathcal{R}$, defined by (1), wherein $\mu$ is taken to be a *completely random measure* without drift on $\mathcal{R}$. The law of $\mu$ is uniquely characterized by the Laplace functional

$$(3) \qquad \mathcal{L}_\mu(g|\rho, \eta) = \exp\left[-\int_\mathcal{R} \int_\mathcal{R} (1 - e^{-g(u)z}) \rho(dz|u) \eta(du)\right],$$

where $g$ is a nonnegative function on $\mathcal{R}$. We say that $\mu$ is *characterized* by an intensity measure $\rho(dz|u)\eta(du)$, as it can be represented in a distributional sense as

$$\mu(du) = \int_\mathcal{R} z \mathcal{N}(dz, du),$$

where $\mathcal{N}(dz, du)$ is a Poisson random measure, taking on points $(z, u)$ in $\mathcal{R} \times \mathcal{R}$, with mean intensity

$$(4) \qquad \mathbb{E}[\mathcal{N}(dz, du)] = \rho(dz|u)\eta(du),$$

such that, for any bounded set $B$ on the half line, $\int_B \int_\mathcal{R} \min(z, 1)\rho(dz|u) \times \eta(du) < \infty$.

Suppose we collect observations $\mathbf{T} = (T_1, \ldots, T_N)$ from $N$ items with hazard rates given by (1) until time $\tau$, so that $T_1 < \cdots < T_n < \tau$ denote the completely observed failure times, and $T_{n+1} = \cdots = T_N \equiv \tau$ are the right-censored times. Define

$$(5) \qquad g_N(u) = \int_0^\tau \left[\sum_{i=1}^N \mathbb{I}(T_i \geq t)\right]\mathbb{I}(t < u)\,dt,$$

and $\mu(g_N) = \int_{\mathcal{R}} g_N(u)\mu(du) = \int_0^\tau [\sum_{i=1}^N \mathbb{I}(T_i \geq t)]\lambda(t|\mu)\,dt$, where $\sum_{i=1}^N \mathbb{I}(T_i \geq t)$ is called the total time transform [4]. Assuming a multiplicative intensity model [1, 2], the likelihood of the data is

$$(6) \qquad \frac{N!}{(N-n)!}\left[\prod_{i=1}^n \int_{\mathcal{R}} \mathbb{I}(T_i < u_i)\mu(du_i)\right]\exp[-\mu(g_N)].$$

(See Remark 2.2.) Noticing that this is a special case of the Lévy moving averages model considered in [32] with a general kernel replaced by the indicator function, we can describe the posterior distribution in terms of partitions of $n$ integers (given in the Appendix). Due to the special structure of the indicator kernel [see (24)], we recognize from the posterior distribution that the information carried by a partition about the remaining members other than the maximal element in any cell is irrelevant. In other words, only the maximal element and the size of each cell is *sufficient* for this problem. To summarize the information, we can define an integer-valued vector $\mathbf{S} = (S_0, S_1, \ldots, S_{n-1}, S_n)$ (see [10, 14]), referred to as an $\mathbf{S}$-path (of $n+1$ coordinates), which satisfies (i) $S_0 = 0$ and $S_n = n$; (ii) $S_j \leq j$, $j = 1, \ldots, n-1$; and (iii) $S_j \leq S_{j+1}$, $j = 1, \ldots, n-1$. A path $\mathbf{S}$ is said to *correspond to* one or many partitions $\mathbf{p}$, provided that (i) labels of the maximal elements of the $n(\mathbf{p})$ cells in $\mathbf{p}$ coincide with locations $j$ at which $S_j > S_{j-1}$, and (ii) the size $e_i$ of the cell $C_i$ with a maximal element $j$ is identical to $m_j = S_j - S_{j-1}$, for all $i = 1, \ldots, n(\mathbf{p})$ (see [10] and [22] for more discussion).

Define $f_N(z, u) = g_N(u)z$. Given the data $\mathbf{T}$, assume that

$$(7) \qquad \kappa_i(e^{-f_N}\rho|u) = \int_{\mathcal{R}} z^i e^{-g_N(u)z}\rho(dz|u) < \infty,$$

for any integer $i \leq n$ and a fixed $u > 0$. Write $\sum_{\mathbf{S}}$ as summing over all paths $\mathbf{S}$ of the same number of coordinates, and $\prod_{\{j\,:\,m_j>0\}}$ and $\sum_{\{j\,:\,m_j>0\}}$ as $\prod_{j=1\,:\,m_j>0}^n$ and $\sum_{j=1\,:\,m_j>0}^n$ conditioned on $\mathbf{S}$, respectively.

THEOREM 2.1.   *Suppose the likelihood of the data is* (6), *and $\mu$ is a completely random measure characterized by the Laplace functional* (3). *Then given the data $\mathbf{T}$, the posterior law of $\mu$ can be described by a three-step experiment:*

(i) *An **S**-path* $\mathbf{S} = (0, S_1, \ldots, S_{n-1}, n)$ *has a distribution* $Z(\mathbf{S}) = \phi(\mathbf{S})/\sum_{\mathbf{S}} \phi(\mathbf{S})$, *where*

$$\phi(\mathbf{S}) = \prod_{\{j\,:\,m_j > 0\}} \binom{j - 1 - S_{j-1}}{j - S_j} \int_{T_j}^{\infty} \kappa_{m_j}(e^{-f_N}\rho|y)\eta(dy). \tag{8}$$

(ii) *Given* $\mathbf{S}$, *there exist* $\sum_{j=1}^{n} \mathbb{I}(S_j > S_{j-1})$ *independent pairs of* $(y_j, Q_j)$, *denoted by* $(\mathbf{y}, \mathbf{Q}) = \{(y_j, Q_j) : m_j > 0, j = 1, \ldots, n\}$, *where* $y_j | \mathbf{S}, \mathbf{T}$ *is distributed as*

$$\eta_j(dy_j|\mathbf{S}, \mathbf{T}) \propto \mathbb{I}(T_j < y_j)\kappa_{m_j}(e^{-f_N}\rho|y_j)\eta(dy_j) \tag{9}$$

*and*

$$\Pr\{Q_j \in dz | \mathbf{S}, y_j, \mathbf{T}\} \propto z^{m_j} e^{-g_N(y_j)z}\rho(dz|y_j). \tag{10}$$

(iii) *Given* $(\mathbf{S}, \mathbf{y}, \mathbf{Q})$, $\mu$ *has the distribution of* $\mu_N^* + \sum_{\{j\,:\,m_j>0\}} Q_j \delta_{y_j}$, *where* $\mu_N^*$ *is a completely random measure characterized by* $e^{-g_N(u)z}\rho(dz|u)\eta(du)$.

REMARK 2.1. The finiteness condition in (7) guarantees the existence of the posterior distributions of $Q_j|\mathbf{S}, y_j$ in statement of (ii) in Theorem 2.1.

COROLLARY 2.1. *Theorem* 2.1 *implies that the posterior mean of the decreasing hazard rate* (1) *given* $\mathbf{T}$ *is given by, for* $t \in [0, \tau]$,

$$\mathbb{E}[\lambda(t|\mu)|\mathbf{T}] = \int_t^{\infty} \kappa_1(e^{-f_N}\rho|y)\eta(dy) + \sum_{\mathbf{S}} Z(\mathbf{S}) \sum_{j=1}^{n} \lambda_j(t|\mathbf{S}), \tag{11}$$

*where* $\kappa_i(e^{-f_N}\rho|y)$, $i = 1, \ldots, n$, *is defined in* (7), $Z(\mathbf{S})$ *is given in Theorem* 2.1 *and*

$$\lambda_j(t|\mathbf{S}) = \frac{\int_{\max(t, T_j)}^{\infty} \kappa_{m_j+1}(e^{-f_N}\rho|y)\eta(dy)}{\int_{T_j}^{\infty} \kappa_{m_j}(e^{-f_N}\rho|y)\eta(dy)} \tag{12}$$

*if* $m_j > 0$, *otherwise* 0.

With the posterior consistency result of Drăgichi and Ramamoorthi [13], the consistency of this Bayes estimate can be obtained via the same argument used in Corollary 1 of [6]. In addition, this path-sum estimator is always less variable than its counterpart in terms of partitions due to the following Rao–Blackwellization result, which states that $\mathbf{p}|\mathbf{S}, \mathbf{T}$ is uniformly distributed over all partitions that correspond to the given path $\mathbf{S}$. (See the proof of Lemma 2.1 in [10] for this total number.)

LEMMA 2.1. *Suppose $\mathbf{S}|\mathbf{T} \sim Z(\mathbf{S})$. Then there exists a conditional distribution*

$$\pi(\mathbf{p}|\mathbf{S},\mathbf{T}) = \frac{1}{\prod_{\{j:m_j>0\}} \binom{j-1-S_{j-1}}{j-S_j}}, \qquad \mathbf{p} \in \mathbb{C}_\mathbf{S}, \tag{13}$$

*where $\mathbb{C}_\mathbf{S}$ is the collection of all partitions that correspond to path $\mathbf{S}$.*

REMARK 2.2. As discussed in [32], Section 4 (see also [3], Section III.2), the multiplicative intensity model captures a large variety of models that appear in event history analysis. Bayesian analysis for models under different censoring schemes, such as left truncation together with right censorship, and random censoring at different time points, follows similarly as the likelihoods differ slightly from (6).

REMARK 2.3. One may model a "bathtub" or $U$-shaped hazard rate with a minimum at $a$ by [41]

$$\lambda(t|a,\mu) = \int \mathbb{I}(|t-a| \geq u)\mu(du),$$

and obtain a posterior distribution characterized by $\mathbf{S}$-paths. In particular, the posterior law of $(a,\mu)$ can then be jointly described by the posterior laws of $a$ and $\mu|a$, where the latter follows naturally as a path characterization for a fixed $a$; see [10] and [23] for estimation in similar mixture models.

**3. The Markov chain Monte Carlo method.** This section introduces an MCMC path-sampler to efficiently compute the hazard estimate (11) (and in general sums over $\mathbf{S}$-paths that appear in many other aforementioned problems). The algorithm samples a Markov chain with a unique stationary distribution $Z(\mathbf{S})$ in a state space as the collection of $\mathbf{S}$-paths of $n+1$ coordinates. It is named an *accelerated path* ($AP$) *sampler* as it is designed to accelerate a straightforward Gibbs sampler [16] in the sense that the algorithm allows more efficient movements among different $\mathbf{S}$-paths.

A straightforward Gibbs sampler, which has a stationary distribution proportional to $\phi(\mathbf{S})$ [see (8)], can be defined [22]: Each Gibbs cycle consists of sampling $S_r|\mathbf{S}_{-r}$, where $\mathbf{S}_{-r} = (S_1,\ldots,S_{r-1},S_{r+1},\ldots,S_{n-1})$ is the "deleted-$r$" vector, and cycling through $r = 1,\ldots,n-1$. The conditional probabilities are

$$\Pr\{S_r = j|\mathbf{S}_{-r}\} \propto \phi(0,S_1,\ldots,S_{r-1},j,S_{r+1},\ldots,S_{n-1},n),$$

for $j = S_{r-1}, S_{r-1}+1,\ldots,S_{r+1}-1, S_{r+1}$ (subject to the definition of an $\mathbf{S}$-path). At step $r$, $S_r$ would remain unchanged if $S_{r-1} = S_{r+1}$. This retards the convergence of the chain to its equilibrium state, and thus results in poor approximations of sums over $\mathbf{S}$-paths. The above phenomenon motivates us to accelerate this naïve chain in accordance with an increasing number of possible movements among the state space within any step. Noticing that each step of the Gibbs sampler is equivalent to re-determinations of the two increments $m_r$ and $m_{r+1}$ at locations $r$ and $r+1$, respectively, our idea is to replace $m_{r+1}$ by some other $m_q$ such that, at any step, it is relatively less likely that the resulting chain is bounded to remain unchanged. Suppose $q > r$ denotes the next location at which $m_q = S_q - S_{q-1} > 0$.

ALGORITHM 3.1 (*Accelerated path sampler*). A Markov chain of **S**-paths (of $n+1$ coordinates), which has a unique stationary distribution $Z(\mathbf{S}) = \phi(\mathbf{S})/\sum_{\mathbf{S}} \phi(\mathbf{S})$ with $\phi(\mathbf{S})$ given by (8), can be defined by a transition cycle of $n-1$ steps:

(i) At step $r$, suppose $\mathbf{S} = (0, S_1, \ldots, S_{r-1}, c, \ldots, c, S_q, \ldots, S_{n-1}, n)$, where $S_{r-1} \leq c \leq \min(r, S_q - 1)$. **S** moves to $(0, S_1, \ldots, S_{r-1}, j, \ldots, j, S_q, \ldots, S_{n-1}, n)$ with conditional probability proportional to

$$
(14) \qquad \frac{r - S_{r-1}}{S_q - 1 - S_{r-1}} \int_{T_q}^{\infty} \kappa_{S_q - S_{r-1}}(e^{-f_N} \rho | y) \eta(dy),
$$

if $j = S_{r-1}$; otherwise, if $j \in \{S_{r-1} + 1, S_{r-1} + 2, \ldots, \min(r, S_q - 1)\}$, with probability proportional to

$$
(15) \qquad \binom{S_q - S_{r-1} - 2}{S_q - j - 1} \prod_{i=r+1}^{q-1} \left( \frac{i - j}{i - S_{r-1}} \right)
$$
$$
\times \int_{T_r}^{\infty} \kappa_{j - S_{r-1}}(e^{-f_N} \rho | y) \eta(dy) \times \int_{T_q}^{\infty} \kappa_{S_q - j}(e^{-f_N} \rho | z) \eta(dz).
$$

(ii) Repeat step (i) for $r = 1, 2, \ldots, n-1$ to complete a cycle.

Starting with an arbitrary path $\mathbf{S}^{(0)}$, and repeating $M$ cycles according to the above scheme, gives a Markov chain $\mathbf{S}^{(0)}, \mathbf{S}^{(1)}, \ldots, \mathbf{S}^{(M)}$ with a unique stationary distribution $Z(\mathbf{S})$. Then, for large $M$, the ergodic average

$$
(16) \qquad \widehat{\lambda}_M(t) = \int_t^{\infty} \kappa_1(e^{-f_N} \rho | y) \eta(dy) + \frac{1}{M} \sum_{i=1}^{M} \sum_{j=1}^{n} \lambda_j(t | \mathbf{S}^{(i)})
$$

approximates the hazard estimate (11) [44].

The validity of the AP sampler is justified by an idea in [20] or [48] (see [22] for a proof in more detail in the gamma process case). One could always define a sequence of reducible transition kernels $\mathbb{P}^{(r)}$, $r = 1, 2, \ldots, d$, that all have the target stationary distribution. Multiplying them in series gives a transition kernel $\mathbb{P} = \mathbb{P}^{(1)} \times \mathbb{P}^{(2)} \times \cdots \times \mathbb{P}^{(d)}$ with the target stationary distribution (from construction). If the chain defined by $\mathbb{P}$ is irreducible, as all states communicate, the target stationary distribution will be unique. At each step $r$, $r = 1, \ldots, n-1 (=d)$, of the AP sampler, the kernel $\mathbb{P}^{(r)}$ is defined by the probability that the chain moves from the path after step $r-1$, $\mathbf{S}_0 = (0, S_1, \ldots, S_{r-1}, c, \ldots, c, S_q, \ldots, S_{n-1}, n)$ to $\mathbf{S}^* = (0, S_1, \ldots, S_{r-1}, j, \ldots, j, S_q, \ldots, S_{n-1}, n)$ for $j \in \{S_{r-1}, S_{r-1} + 1, \ldots, \min(r, S_q - 1)\}$, which is proportional to $\phi(\mathbf{S}^*)$. That is, $\mathbf{S}_0$ communicates only with paths in the collection defined by $\{\mathbf{S} : \mathbf{S} = (0, S_1, \ldots, S_{r-1}, j, \ldots, j, S_q, \ldots, S_{n-1}, n), j \in \{S_{r-1}, S_{r-1}+1, \ldots, \min(r, S_q-1)\}\}$. With this construction, the kernel $\mathbb{P}^{(r)}$ decomposes the state space of all $\mathbf{S}$-paths of $n+1$ coordinates into a finite collection of mutually exclusive communication classes (see Theorem 3 in [15], page 392). One can easily check that each kernel $\mathbb{P}^{(r)}$, though not irreducible, has a stationary distribution $Z(\mathbf{S})$. More importantly, the chain defined by $\mathbb{P} = \mathbb{P}^{(1)} \times \mathbb{P}^{(2)} \times \cdots \times \mathbb{P}^{(n-1)}$ is irreducible, as all states can communicate with the path $\mathbf{S} = (0, 0, \ldots, 0, n)$ within one cycle. Hence, the AP sampler gives a Markov chain of $\mathbf{S}$-paths with a unique stationary distribution $Z(\mathbf{S})$.

**4. Examples.** One can model $\mu$ in (1) by a variety of random measures. Corresponding posterior analyses follow from the results in the previous sections. This section looks at two explicit examples wherein $\mu$ is characterized by the mean measure,

$$(17) \quad \rho_{\alpha,\beta}(dz|u)\eta(du) = \frac{1}{\Gamma(1-\alpha)} z^{-\alpha-1} \exp[-z/\beta(u)] \, dz \, \eta(du).$$

This class of random measures generalizes the generalized gamma random measure proposed by Brix [8], for $0 < \alpha < 1$ and $0 \le \beta < \infty$, or, $-\infty < \alpha \le 0$ and $0 < \beta < \infty$. It includes the weighted gamma process (when $\alpha = 0$), a stable law with index $0 < \alpha < 1$ (when $\beta = \infty$), and the inverse-Gaussian process (when $\alpha = 1/2$ and $\beta > 0$) (see [32]).

4.1. *The weighted gamma random measure.* If $\alpha = 0$ in (17), $\mu$ is the weighted gamma random measure with shape measure $\eta$ and scale measure $\beta$. Corollary 2.1 gives the Bayes estimate of the decreasing hazard rate (1) according to

$$(18) \quad \kappa_i(e^{-f_N}\rho|u) = (i-1)! \times [\beta(u)^{-1} + g_N(u)]^{-i}, \qquad i = 1, 2, \ldots, n.$$

To apply the AP sampler, one needs to compute conditional probabilities (14) and (15) that are proportional to

$$(r - S_{r-1}) \times \xi_{S_q - S_{r-1}}(T_q|\mathbf{T})$$

and

$$\left[\prod_{i=r+1}^{q-1} \left(\frac{i-j}{i-S_{r-1}}\right)\right] \times \xi_{j-S_{r-1}}(T_r|\mathbf{T}) \times \xi_{S_q-j}(T_q|\mathbf{T}),$$

respectively, where $\xi_i(t|\mathbf{T}) = \int_t^\infty [\beta(v)^{-1} + g_N(v)]^{-i} \eta(dv)$, $i = 1, \ldots, n$.

4.2. *The stable law.* The stable law with index $0 < \alpha < 1$ appears when $\beta = \infty$ in (17). The posterior mean of (1) defined by this class of random measures follows from Corollary 2.1, and it can be evaluated by implementing the AP sampler, based on

$$\kappa_i(e^{-f_N}\rho|y) = \frac{\Gamma(i-\alpha)}{\Gamma(1-\alpha)[g_N(y)]^{i-\alpha}}, \qquad i=1,2,\ldots,n.$$

**5. Numerical results.** This section addresses the effectiveness of the AP sampler for evaluating the hazard estimate (11). In particular, we select a special case of (1) wherein $\mu$ is a gamma process. The posterior analysis follows from discussions in Section 4.1 with $\beta(\cdot) = 1$. Hence, the posterior mean of the monotone hazard rate reduces to

$$(19) \quad \mathbb{E}[\lambda(t|\mu)|\mathbf{T}] = \xi_1(t|\mathbf{T}) + \sum_{\mathbf{S}} Z^*(\mathbf{S}) \sum_{j=1}^n \left[ m_j \frac{\xi_{m_j+1}(\max(t,T_j)|\mathbf{T})}{\xi_{m_j}(T_j|\mathbf{T})} \right],$$

where $Z^*(\mathbf{S}) \propto \prod_{\{j\,:\,m_j>0\}} (j-1-S_{j-1})!/(j-S_j)! \times \xi_{m_j}(T_j|\mathbf{T})$. The complexity in evaluating $\xi_j(t|\mathbf{T})$ can be reduced by assuming a uniform shape probability from 0 to $6(\geq \tau)$ for $\eta(\cdot)$, even though the closed-form expression is tedious (see [25] for its exact expression). The methodology is tested by data from a piecewise constant hazard rate model, for which the hazard rate of an item is

$$(20) \qquad \lambda(t) = \begin{cases} 1, & 0 \leq t < 1, \\ 0.5, & t \geq 1. \end{cases}$$

Data are generated subject to a termination time $\tau = 3$, such that the censoring rate is about 15%. All simulation results that follow are based on a Monte Carlo size $M = 1000$, and an initial path $\mathbf{S}^{(0)} = (0,1,\ldots,n-1,n)$.

REMARK 5.1. In practice, the implementation of the AP sampler depends heavily on evaluations of the double integral

$$\int_y^\infty \int_\mathcal{R} z^i e^{-g_N(u)z} \rho(dz|u)\,\eta(du), \qquad y>0, i=1,2,\ldots,n,$$

where, according to (5),

$$g_N(u) = \begin{cases} \sum_{i=1}^{j-1} T_i + (N-j+1)u, & T_{j-1} < u \leq T_j, j=1,\ldots,n+1, \\ \sum_{i=1}^n T_i + (N-n)\tau, & u > \tau. \end{cases}$$

It is important to note that the inner integral, which is defined to be $\kappa_i(e^{-f_N}\rho|u)$ in (7), is the conditional cumulant of an (exponentially tilted) infinitely divisible random variable with density of an available form for any given $\rho$. Then it follows that the integral may be calculated using a result of T. N. Thiele, which gives a recursive relation between cumulants and moments of a random variable (see, e.g., [18] and [43], Section 2.3). See [33], Section 4.1 and [34], Section 5, for more discussion of this problem appearing in other contexts.

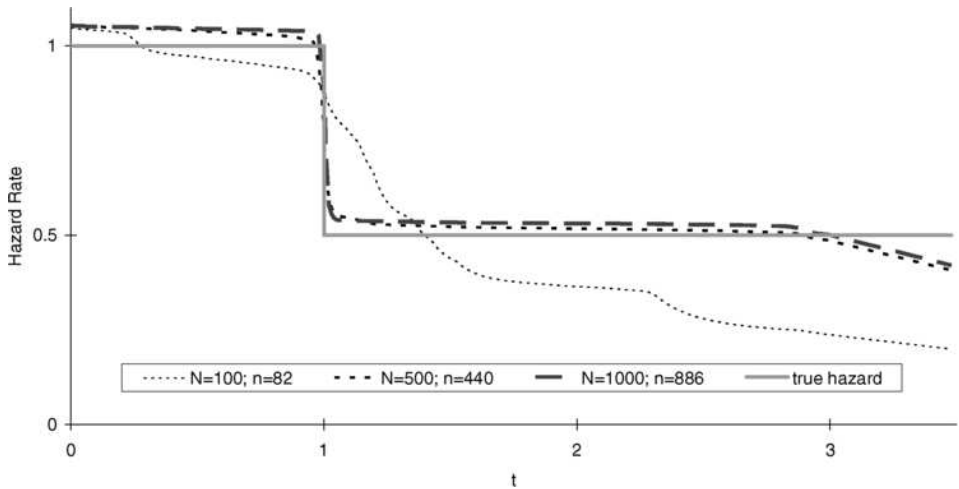

Fig. 1. *MCMC hazard estimates produced by the AP sampler.*

5.1. *Resolution of the AP sampler.* This section focuses on the convergence property of the hazard estimate (19) approximated by the AP sampler as the sample size $N$ increases. Based on nested samples of sizes $N = 100$, 500 and 1000, MCMC estimates (16) according to (18) are displayed in Figure 1. The graphs echo the fact that the approximated Bayes estimate of the decreasing hazard function, $\widehat{\lambda}_M(t)$, tends to the "true" hazard rate (20) as sample size increases. We remark that the drop of the hazard estimates after $t = 3$ results from the fact that the estimates are mainly constructed based on the prior information, as no complete data is observed after that time point.

5.2. *Comparison with other methods.* Path-sum estimates of monotone hazard rates, though they appeared two decades ago, have not received much attention and are not commonly used due to unavailability of efficient numerical methods; previous attempts by Brunner and Lo [10, 11] and Brunner [9] are all far from successful due to their incapability of sampling from desired posterior distributions of **S**-paths in the respective models. On the contrary, partition-sum counterparts have been used as a substitute [21, 25] since there are many well-developed numerical methods for sampling partitions (see, e.g., [28, 29, 30, 40]). This section aims at comparing MCMC hazard estimates produced by the AP sampler, the Gibbs path sampler defined in [22] and a Gibbs sampler for partitions defined in [40] (see [25] for an exact description of the algorithm being applied to this gamma model). Standard errors of MCMC hazard estimates by the three different methods are estimated by repetitions of experiment, and are used as the standard of comparison.

TABLE 1
*A large-sample study of MCMC hazard estimates according to* 1000 *independent replications of the accelerated path (AP) sampler, of the Gibbs path (gP) sampler and of the weighted Chinese restaurant process (gWCR) sampler*

| $t$ | MCMC method | Average of hazard estimates at time $t$ | Standard error of hazard estimates at time $t$ |
|---|---|---|---|
|  | AP | 0.9667340 | 0.0038426 |
| 0.5 | gP | 0.9677822 | 0.0517161 |
|  | gWCR | 0.9668593 | 0.0080103 |
|  | AP | 0.8815065 | 0.0065156 |
| 0.99 | gP | 0.8812967 | 0.0594452 |
|  | gWCR | 0.8820966 | 0.0097827 |
|  | AP | 0.8530503 | 0.0067767 |
| 1.01 | gP | 0.8524295 | 0.0541091 |
|  | gWCR | 0.8537771 | 0.0099268 |
|  | AP | 0.3708132 | 0.0055500 |
| 2.0 | gP | 0.3692810 | 0.0440281 |
|  | gWCR | 0.3707660 | 0.0106327 |

Here the sample size is fixed at $N = 100$, and there are $n = 82$ complete observations in our simulated data set. Markov samples from all the three Markov chain experiments are collected after a "burn-in" period of 10000 cycles. We compute 1000 independent hazard estimates by 1000 repetitions of each experiment. These are used to estimate the average and the standard error of the hazard estimates in the usual manner. The hazard rates $\lambda(t)$ at times $t = 0.5$, 0.99, 1.01 and 2.0 are studied. The points, 0.99 and 1.01, are near the change point at 1, and they seem to reflect well the effectiveness and the efficiency of the MCMC methods in the worst case. Table 1 displays the averages and the standard errors of the 1000 realizations of the MCMC hazard estimates produced by the three methods. At different time points, the three averages are close to each other, yet the standard errors vary substantially. The standard error of hazard estimates produced by the AP sampler is the smallest among all the three methods. On one hand, the AP sampler definitely outweighs the naïve Gibbs path sampler. On the other hand, the AP sampler beats the closest competitor, the gWCR sampler, by a comfortable margin. These show that our "acceleration" scheme works extremely well.

**6. Proportional hazards.** The Cox regression model [12] is an important example of the multiplicative intensity model that can allow incorporation of covariates, together with right independent censoring, in survival analysis. For Bayes inference of general hazard rates with the presence of covariates, see [27, 30, 31, 35], among others. Assume that the underlying hazard defined on $\mathcal{R}$ is modeled by

$$\lambda(t|\mathbf{Z},\boldsymbol{\theta},\mu) = \int_{\mathcal{R}} \exp(\boldsymbol{\theta}^T\mathbf{Z})\mathbb{I}(t<u)\mu(du),$$

where $\mathbf{Z}$ is a covariate vector with parameter vector $\boldsymbol{\theta}$, and $\lambda_0(t|\mu) = \int_{\mathcal{R}} \mathbb{I}(t<u)\mu(du)$, same as (1), is a decreasing baseline hazard rate. Suppose we collect data until time $\tau$ and the data $\mathbf{D} = ((T_1, \mathbf{Z}_1), \ldots, (T_N, \mathbf{Z}_N))$ summarize completely observed failure times $T_1 < \cdots < T_n$ and right-censored times $T_i = \tau$, $i = n+1, \ldots, N$, associated with covariate vectors $\mathbf{Z}_i$, $i = 1, \ldots, N$, with unknown parameter vector $\boldsymbol{\theta}$. Define $f_{N,\boldsymbol{\theta}}(z,u) = g_{N,\boldsymbol{\theta}}(u)z$, where

$$(21) \qquad g_{N,\boldsymbol{\theta}}(u) = \int_0^\tau \left[\sum_{i=1}^N \mathbb{I}(T_i \geq t)\exp(\boldsymbol{\theta}^T\mathbf{Z}_i)\right]\mathbb{I}(t<u)\, dt.$$

Then the Cox proportional hazards likelihood may be written as

$$(22) \qquad \left[\prod_{i=1}^n \exp(\boldsymbol{\theta}^T\mathbf{Z}_i)\lambda_0(T_i|\mu)\right]\exp[-\mu(g_{N,\boldsymbol{\theta}})],$$

where $\mu(g_{N,\boldsymbol{\theta}}) = \int_{\mathcal{R}} g_{N,\boldsymbol{\theta}}(u)\mu(du) = \int_0^\tau [\sum_{i=1}^N \mathbb{I}(T_i \geq t)\exp(\boldsymbol{\theta}^T\mathbf{Z}_i)]\lambda_0(t|\mu)\, dt$. Assume $\int_{\mathcal{R}} z^i e^{-g_{N,\boldsymbol{\theta}}(u)z}\rho(dz|u) < \infty$, for $i = 1, \ldots, n$ and a fixed $u > 0$.

PROPOSITION 6.1. *Suppose the likelihood of the data is given by* (22). *Let $\pi(d\boldsymbol{\theta})$ denote our prior for $\boldsymbol{\theta}$ and independently assume $\mu$ is a completely random measure characterized by the Laplace functional* (3). *Then the posterior distribution of $\mu|\boldsymbol{\theta},\mathbf{D}$ can be described as a three-step hierarchical experiment as in Theorem* 2.1, *of which $f_N(\cdot,\cdot)$ and $g_N(\cdot)$ are replaced by $f_{N,\boldsymbol{\theta}}(\cdot,\cdot)$ and $g_{N,\boldsymbol{\theta}}(\cdot)$, respectively.*

To evaluate any posterior quantities of model (22), such as the posterior mean of the underlying monotone baseline hazard rate and the posterior mean of the covariate parameters $\boldsymbol{\theta}$, run the following Gibbs sampler:

1. Draw $\mathbf{S}|\mathbf{Q},\mathbf{y},\boldsymbol{\theta},\mathbf{D}$ by implementing Algorithm 3.1 with $f_N(\cdot,\cdot)$ and $g_N(\cdot)$ replaced by $f_{N,\boldsymbol{\theta}}(\cdot,\cdot)$ and $g_{N,\boldsymbol{\theta}}(\cdot)$, respectively.

2. Draw $\mathbf{Q}, \mathbf{y} | \mathbf{S}, \boldsymbol{\theta}, \mathbf{D}$ according to the analogues of the conditional distributions (9) and (10) in Theorem 2.1 with $f_N(\cdot, \cdot)$ and $g_N(\cdot)$ replaced by $f_{N,\boldsymbol{\theta}}(\cdot, \cdot)$ and $g_{N,\boldsymbol{\theta}}(\cdot)$, respectively.
3. Draw $\boldsymbol{\theta} | \mathbf{Q}, \mathbf{y}, \mathbf{S}, \mathbf{D}$ from the density proportional to

$$\pi(d\boldsymbol{\theta})B(\boldsymbol{\theta}) \prod_{\{j\,:\,m_j>0\}} e^{-g_{N,\boldsymbol{\theta}}(y_j)Q_j},$$

where

(23) $$B(\boldsymbol{\theta}) = \exp\left[-\int_\mathcal{R}\int_\mathcal{R}(1-e^{-g_{N,\boldsymbol{\theta}}(u)z})\rho(dz|u)\eta(du)\right]\prod_{i=1}^n \exp(\boldsymbol{\theta}^T\mathbf{Z}_j).$$

REMARK 6.1. Note that $g_{N,\boldsymbol{\theta}}(u)$ is again a piecewise linear function of $u$ as $g_N(u)$ in the case without covariates (discussed in Remark 5.1). This does not create any further complexities in evaluating integrals at steps 1 and 2 of the Gibbs sampler. Step 3, which is of the same form as the step 4 (for conditional draws of regression parameters $\boldsymbol{\theta}$) of the blocked Gibbs algorithm suggested by Ishwaran and James ([30], page 184), can be dealt with via a Metropolis step.

REMARK 6.2. We conclude here that $\mathbf{S}$-paths may be derived from every exchangeable random partition $\mathbf{p}$ by summation. That is to say, the general correspondence between these structures is simply a combinatorial relationship. However, what we are exploiting is the fact that the monotone hazard rates models, due to the presence of the indicator kernel, are naturally representable in terms of tractable $\mathbf{S}$-path structures. Thus, technically our approach may be applied to models exhibiting similar structure (see, e.g., [10, 11, 23, 24]). With this point in mind, we note that recently James and Lau [33] show that the non-Gaussian Ornstein–Uhlenbeck models of Barndorff-Nielsen and Shephard [5], which are also special types of Lévy moving averages, also exhibit $\mathbf{S}$-path structures which are amenable to our approach. It is important to note that those models are not of the form in (1).

## APPENDIX

PROOF OF THEOREM 2.1 AND LEMMA 2.1. The proof relies on the following two key structures:

($\star$) For any $\mathbf{p} \in \mathbb{C}_\mathbf{S}$,

(24) $$\prod_{i=1}^{n(\mathbf{p})} \mathbb{I}\left(\max_{j \in C_i} T_j < v_i\right) = \prod_{\{j\,:\,m_j>0\}} \mathbb{I}(T_j < y_j),$$

where $\mathbf{v} = (v_1, \ldots, v_{n(\mathbf{p})})$ represents the unique values among $(u_1, \ldots, u_n)$ in (6) and $\mathbf{y} = \{y_j : m_j > 0, j = 1, \ldots, n\}$ is a permutation of $\mathbf{v}$ according to $\mathbf{p} \in \mathbb{C}_\mathbf{S}$.

($\star\star$) The total number of partitions that correspond to a given $\mathbf{S}$ equals

$$(25) \qquad |\mathbb{C}_\mathbf{S}| := \sum_{\mathbf{p} \in \mathbb{C}_\mathbf{S}} 1 = \prod_{\{j : m_j > 0\}} \binom{j - 1 - S_{j-1}}{j - S_j},$$

where $\sum_{\mathbf{p} \in \mathbb{C}_\mathbf{S}}$ represents summing over all partitions that correspond to $\mathbf{S}$.

According to James [32], the posterior law of $\mu | \mathbf{T}$ is equivalent to the distribution of a random measure $\mu_N^* + \sum_{i=1}^{n(\mathbf{p})} J_i \delta_{v_i}$. It is determined by the joint distribution of $\mu_N^*, \mathbf{J}, \mathbf{v}, \mathbf{p} | \mathbf{T}$, which is proportional to

$$(26) \qquad \mathbb{P}(d\mu_N^*) \prod_{i=1}^{n(\mathbf{p})} J_i^{e_i} e^{-g_N(v_i) J_i} \rho(dJ_i | v_i) \prod_{i=1}^{n(\mathbf{p})} \mathbb{I}\left(\max_{j \in C_i} T_j < v_i\right) \eta(dv_i),$$

where $\mathbf{J} = (J_1, \ldots, J_{n(\mathbf{p})})$, and $\mathbb{P}(d\mu_N^*)$ is a completely random measure characterized by an intensity measure $e^{-g_N(u) z} \rho(dz | u) \eta(du)$ as in (iii) of Theorem 2.1.

Notice that (26), due to its irrelevance to the remaining members other than the maximal elements of the cells in a partition, may be rewritten in terms of the intrinsic characteristics of a path $\mathbf{S}$, provided that $\mathbf{p} \in \mathbb{C}_\mathbf{S}$ based on ($\star$) as

$$(27) \quad \mathbb{P}(d\mu_N^*) \prod_{\{j : m_j > 0\}} Q_j^{m_j} e^{-g_N(y_j) Q_j} \rho(dQ_j | y_j) \prod_{\{j : m_j > 0\}} \mathbb{I}(T_j < y_j) \eta(dy_j),$$

where $\mathbf{Q} = \{Q_j : m_j > 0, j = 1, \ldots, n\}$, is a relabeling of $\mathbf{J}$ according to the correspondence $\mathbf{p} \in \mathbb{C}_\mathbf{S}$. That is, the conditional law of $\mu_N^*, \mathbf{J}, \mathbf{v}, \mathbf{p} | \mathbf{T}$ only depends on **p** through **S**. Equation (27) and the equivalence in distribution relation between the two random measures,

$$\mathcal{L}\left\{\mu_N^* + \sum_{i=1}^{n(\mathbf{p})} J_i \delta_{v_i} \middle| \mathbf{T}\right\} \stackrel{d}{=} \mathcal{L}\left\{\mu_N^* + \sum_{\{j\,:\,m_j>0\}} Q_j \delta_{y_j} \middle| \mathbf{T}\right\}, \quad (28)$$

imply that the posterior law of $\mu|\mathbf{T}$ can be expressed in terms of $\mu_N^*|\mathbf{Q}, \mathbf{y}, \mathbf{S}, \mathbf{T}$ and $\mathbf{Q}, \mathbf{y}, \mathbf{S}|\mathbf{T}$. Then, integrating out $\mu_N^*$ and summing over all $\mathbf{p} \in \mathbb{C}_{\mathbf{S}}$ in (27) yields that the distribution of $\mathbf{Q}, \mathbf{y}, \mathbf{S}|\mathbf{T}$ is proportional to

$$\prod_{\{j\,:\,m_j>0\}} Q_j^{m_j} e^{-g_N(y_j)Q_j} \rho(dQ_j|y_j) \prod_{\{j\,:\,m_j>0\}} \mathbb{I}(T_j < y_j)\eta(dy_j) \left[\sum_{\mathbf{p} \in \mathbb{C}_{\mathbf{S}}} 1\right], \quad (29)$$

where $\sum_{\mathbf{p} \in \mathbb{C}_{\mathbf{S}}} 1$ is given in $(\star\star)$. Now, the laws of $\mathbf{Q}|\mathbf{S}, \mathbf{y}, \mathbf{T}$, $\mathbf{y}|\mathbf{S}, \mathbf{T}$ and $\mathbf{S}|\mathbf{T}$ follow from Bayes' theorem and multiplication rule. Hence, the result in Theorem 2.1 follows, while the conditional distribution of $\mathbf{p}|\mathbf{S}, \mathbf{T}$ in Lemma 2.1 follows by dividing (26) without the leading term $\mathbb{P}(d\mu_N^*)$ by (29). □

TABLE 2
*Total numbers of $\mathbf{S}$-paths and partitions versus sample sizes $(n)$*

| $n$ | # of paths | # of partitions | Ratio in % |
|---|---|---|---|
| 1 | 1 | 1 | 100.000 |
| 3 | 5 | 5 | 100.000 |
| 5 | 42 | 52 | 80.769 |
| 7 | 429 | 877 | 48.917 |
| 10 | 16,796 | 115,975 | 14.482 |
| 15 | 9,694,845 | 1,382,958,545 | 0.701 |
| 20 | 6,564,120,420 | 51,724,158,235,372 | 0.013 |

PROOF OF PROPOSITION 6.1. Following the same arguments as in [32] in getting the posterior distribution given by (26) yields that the posterior law of $\mu|\mathbf{D}$ is equivalent to the distribution of a random measure $\mu_N^* + \sum_{i=1}^{n(\mathbf{p})} J_i \delta_{v_i}$. It is determined by the joint distribution of $\mu_N^*, \mathbf{J}, \mathbf{v}, \mathbf{p}, \boldsymbol{\theta}|\mathbf{D}$, which is proportional to

$$\mathbb{P}(d\mu_N^*)\pi(d\boldsymbol{\theta})B(\boldsymbol{\theta}) \prod_{i=1}^{n(\mathbf{p})} J_i^{e_i} e^{-g_{N,\boldsymbol{\theta}}(v_i)J_i} \rho(dJ_i|v_i)$$

$$\times \prod_{i=1}^{n(\mathbf{p})} \mathbb{I}\left(\max_{j \in C_i} T_j < v_i\right)\eta(dv_i), \quad (30)$$

where $\mathbb{P}(d\mu_N^*)$ is a completely random measure characterized by an intensity measure $e^{-g_{N,\boldsymbol{\theta}}(u)z}\rho(dz|u)\eta(du)$, and $B(\boldsymbol{\theta})$ is given by (23). As in the proof of Theorem 2.1, $\mu_N^*|\mathbf{Q}, \mathbf{y}, \mathbf{S}, \boldsymbol{\theta}, \mathbf{D}$ is equivalent to $\mathbb{P}(d\mu_N^*)$. Then, summing over all $\mathbf{p} \in \mathbb{C}_{\mathbf{S}}$ in (30) yields that the conditional distribution of $\mathbf{Q}, \mathbf{y}, \mathbf{S}, \boldsymbol{\theta}|\mathbf{D}$ is proportional to

$$\pi(d\boldsymbol{\theta})B(\boldsymbol{\theta}) \prod_{\{j\,:\,m_j>0\}} Q_j^{m_j} e^{-g_{N,\boldsymbol{\theta}}(y_j)Q_j} \rho(dQ_j|y_j)$$

$$\times \prod_{\{j\,:\,m_j>0\}} \mathbb{I}(T_j < y_j)\eta(dy_j) \times |\mathbb{C}_{\mathbf{S}}|, \quad (31)$$

due to $(\star)$ and $(\star\star)$. Hence, by Bayes' theorem and multiplication rule, the result follows from the conditional distribution of $\mathbf{Q}, \mathbf{y}, \mathbf{S}|\boldsymbol{\theta}, \mathbf{D}$, which is proportional to (31) without the leading term $\pi(d\boldsymbol{\theta})B(\boldsymbol{\theta})$. □

**Acknowledgments.** This work, partially done when the author was at the Hong Kong University of Science and Technology, was an extension of the author's Ph.D. thesis written under the supervision of Albert Lo. The idea to study this extension was suggested by Lancelot James. The author has greatly benefited from discussions with them in preparing the manuscript. Thanks are also due to the referees for their helpful suggestions.

Department of Statistics
and Applied Probability
National University of Singapore
6 Science Drive 2
Singapore 117546
Republic of Singapore
E-mail: stahmw@nus.edu.sg